\documentclass[12pt,reqno,a4paper]{amsart}
\usepackage[utf8]{inputenc}
\usepackage[T1]{fontenc}

\usepackage{verbatim}

\usepackage[mathscr]{eucal}

\usepackage{amsmath}
\usepackage{amssymb}
\usepackage{amsthm}
\usepackage{lmodern}
\usepackage{latexsym}
\usepackage{amsfonts}
\usepackage{mathrsfs}
\usepackage[all]{xy}
\usepackage{bm}
\usepackage{yfonts}

\usepackage{fge}
\usepackage{hyperref}
\hypersetup{pdftoolbar=true, pdftitle={Proj_Ellip}, pdffitwindow=true, colorlinks=true, citecolor=blue, filecolor=black, linkcolor=blue, urlcolor=blue, hypertexnames=false}

\numberwithin{equation}{section}

\usepackage[top=4cm, bottom=3.4cm, left=3cm , right=3cm]{geometry}

\allowdisplaybreaks[4]

\usepackage{enumerate}

\theoremstyle{plain}
\newtheorem{thm}{Theorem}[section]
\newtheorem{lemma}[thm]{Lemma}
\newtheorem{coro}[thm]{Corollary}

\newcommand{\bee}{\begin{eqnarray}}
\newcommand{\ene}{\end{eqnarray}}
\newcommand{\bea}{\begin{eqnarray*}}
\newcommand{\ena}{\end{eqnarray*}}

\usepackage{fancyhdr}
\DeclareFontFamily{U}{mathx}{}
\DeclareFontShape{U}{mathx}{m}{n}{ <-> mathx10 }{}
\DeclareSymbolFont{mathx}{U}{mathx}{m}{n}
\DeclareFontSubstitution{U}{mathx}{m}{n}

\DeclareMathAccent{\widecheck}{0}{mathx}{"71}

\title{simultaneous non-vanishing of central values of $\mathrm{GL}_2\times \mathrm{GL}_3$ and $\mathrm{GL}_3\times \mathrm{GL}_3$ $L$-functions}



\author[J-J. PAN]{JUNJIE PAN}
\address{School of Sciences, Xi'an University of Technology,
Xi'an 710054, China} \email{junjiepan24@163.com}

\begin{document}

\maketitle

\begin{abstract}
\small{Let $g$ denote a fixed holomorphic Hecke cusp form of weight $k \equiv 0 \pmod{4}$ on $\mathrm{SL}_2(\mathbb{Z})$, 
and let $\pi$ be a fixed cuspidal automorphic representation of $\mathrm{GL}_3$. 
In this paper, we establish an asymptotic formula for the first moment of the product
\[L(1/2,g\times F)L(1/2,\pi\times F),\]
where $F$ runs over an orthonormal basis of Hecke-Maa{\ss} cusp forms 
of level $q$ on $\mathrm{GL}_3$.
As an application, we deduce that $L(1/2, g\times F)L(1/2,\pi\times F)\neq 0$
for infinitely many such forms $F$.
}
\end{abstract}

\section{Introduction}\label{secIntrod}

\setcounter{equation}{0}
\medskip

The behavior of automorphic $L$-functions at special points ranks among the central topics in modern analytic number theory.
In particular, non-vanishing of central values has important arithmetic consequences, including applications to the Gan-Gross-Prasad conjectures and bounds toward the generalized Ramanujan conjecture \cite{LRS}. 
Moreover, non-vanishing results reveal subtle interactions between different automorphic families.

One of the earliest results in this direction is due to Ramakrishnan and Rogawski \cite{RR}, 
who employed the relative trace formula to establish a simultaneous non-vanishing result for the central 
values of the $\mathrm{GL}_2\times \mathrm{GL}_1$ and $\mathrm{GL}_2$ $L$-functions. They proved that
\[
\sum_{\substack{\phi \in \mathcal{F}_N(k) \\ a_p(\phi)\in J}}\frac{1}{(\phi,\phi)}L(1/2,\phi \times \chi)L(1/2,\phi)=2\mu_p(J)c_k L(1,\chi)+O\left(N^{-k/2+\varepsilon}\right),
\]
where $\phi$ runs over newforms of weight $k$ for $\Gamma_0(N)$, and $\chi$ denotes a Dirichlet character.
Building on this approach, Xu \cite{Xu} and Liu \cite{Liu} extended the method to products of two $\mathrm{GL}_2$ $L$-functions,
treating the Maa{\ss} form and holomorphic cusp form, respectively.

The first higher-rank result is due to Li \cite{Li},
who investigated $\mathrm{GL}_3\times \mathrm{GL}_2$ and $\mathrm{GL}_2$ $L$-functions and proved that
\[
\sideset{}{'}\sum_{j}e^{-\frac{t_j^2}{T^2}}L(1/2,F \times u_j)L(1/2,u_j)=\frac{12L(1,F)L(1,\overline{F})}{\pi^3}T^2+O_{\varepsilon,f}\left(T^{11/6+\varepsilon}\right),
\]
which in particular implies that $L(1/2,F\times u_j)L(1/2,u_j)\neq 0$ for infinitely many $u_j$. 
Soon after, Khan \cite{Kh} obtained an analogous result for twists by holomorphic cusp forms of prime level, showing that
\begin{align*}
&\frac{12}{q(k-1)}\sum_{f\in H_k^*(p)}L(1/2,g \times f)L(1/2,f)=\\
&\hskip3cm\frac{L(1,g)L(1,\tilde{g})}{\zeta(2)}\left(1+O_g(k^{-1})\right)+O\left(q^{\frac{2\theta_1-1}{2}+\varepsilon}+q^{\frac{2\theta_2-1}{4}+\varepsilon}+q^{-\frac{1}{8}+\varepsilon}\right),
\end{align*}
where $\theta_1,\theta_2$ are bounds toward the Ramanujan conjecture for $\mathrm{GL}_3$. 
More recently, Maiti and Mallesham \cite{MM} considered the first moment of $L(1/2,g \times F)L(1/2,F)$, 
where $F$ runs over an orthonormal basis of $\mathrm{GL}_3$ Hecke-Maa{\ss} cusp forms of level $q$ on $\mathrm{GL}_3$. They proved that
\begin{align*}
\sum_{F}L(1/2,g \times F)L(1/2,F)=
\frac{1}{195\pi^5}\iint_{\mathrm{Re}(\mu)}M(\mu,k)h(\mu)spec(\mu)d\mu+O\left(T^{\frac{17}{6}+\varepsilon}R^2\right),
\end{align*}
where the main term is $\asymp T^3R^2$, and $g$ is a holomorphic Hecke cusp form for $\mathrm{SL}_2(\mathbb{Z})$.
Regarding the simultaneous non-vanishing of two distinct Rankin-Selberg $L$-functions, 
Blomer and Mili\'cevi\'c \cite{BM} obtained, for two fixed Hecke-Maa{\ss} forms $f_1,f_2$ satisfying the Ramanujan-Petersson conjecture,
\begin{align*}
&\sideset{}{^*}\sum_{\chi \bmod q}L(1/2,f_1 \times \chi)\overline{L(1/2,f_2 \times \chi)}=\\
&\hskip4cm\frac{2}{\zeta(2)}\psi(q)M(f_1,f_2,q)+O_{f_1,f_2}\left(q^{1+\varepsilon}(\mathcal{Q}_1^{-\frac{1}{22}}+(q/\mathcal{Q}_1^2)^{-\frac{1}{22}})\right),
\end{align*}
which implies the existence of characters $\chi$ such that both central values are non-zero.
This refined an earlier result of Hoffstein and Lee \cite{HL} that could not specify the modulus of the character.

In this paper, we consider the first moment of the product of $\mathrm{GL}_2\times \mathrm{GL}_3$ and $\mathrm{GL}_3\times \mathrm{GL}_3$ $L$-functions.
Specifically, we fix a holomorphic Hecke cusp form $g$ for $\mathrm{SL}_2(\mathbb{Z})$, and a Hecke-Maass cusp form $\pi$ for $\mathrm{SL}_3(\mathbb{Z})$.
We study the first moment of $L(s,g\times F)L(s, \pi\times F)$ at the central point as $F$ runs over an orthonormal basis of the space of Hecke-Maa{\ss} cusp forms for $\mathrm{SL}_3(\mathbb{Z})$.

Let $\varpi$ be an $L^2$-normalized newvector in an everywhere spherical cuspidal automorphic representation of
$F\subseteq L^2_{\mathrm{cusp}}(\Gamma_0(q)\backslash\mathbb{H}^3)$,
which is a $\mathrm{GL}_3$-automorphic form of level $q$. 
Let $\mu=\mu_F=(\mu_1,\mu_2,\mu_3)$ denote the Langlands parameters of $F$,
normalized by $A_{\varpi}(1,1)=1$. Under the Ramanujan conjecture, we have $\mu\in(i\mathbb{R})^3$.
We denote by $\mathfrak{a}_{\mathbb{C}}^{*}$ the complexified dual of the Lie algebra of the maximal torus in $\mathrm{PGL}_3(\mathbb{Z})$.
Let $f$ be a normalized Hecke-Maa{\ss} cusp form of type $(v_1,v_2)\in\mathbb{C}^3$, its Fourier-Whittaker expansion takes the form
\begin{align*}
f(z)=\sum_{\gamma\in\mathrm{U}_2(\mathbb{Z})\backslash\mathrm{SL}_2(\mathbb{Z})}\sum_{n_1=1}^{\infty}\sum_{n_2\neq 0}\frac{A_{f}(n_1,n_2)}{|n_1n_2|}W_{J}
\left(\begin{pmatrix}
|n_1n_2|& & \\ & n_1 & \\ & & 1
\end{pmatrix}
\begin{pmatrix}
\gamma &\\ & 1
\end{pmatrix}z,v,\psi_{1,\frac{n_2}{|n_2|}}\right),
\end{align*}
where $\mathrm{U}_{2}(\mathbb{Z})$ denotes the group of upper triangular unipotent matrices in $\mathrm{SL}_2(\mathbb{Z})$,
and $W_{J}$ denotes the Jacquet Whittaker function and 
$\psi_{v_1,v_2}\left(\left(
\begin{smallmatrix}
1&u_2&u_{1,3}\\ & 1 & u_1 \\ & & 1
\end{smallmatrix}
\right)\right)
=e^{2\pi i(v_1u_2+v_2u_1)}$; see \cite[Section 6.2]{Go}.
Let $\pi \subseteq L^2_{\mathrm{cusp}}(\Gamma_0(q)\backslash\mathbb{H}^3)$ be the cuspidal automorphic representation, generated by the Hecke-Maa{\ss} forms $f$,
and set $A_{\pi}(m,n)=A_{f}(m,n)$ for $m,n>0$.
The Rankin-Selberg $L$-function attached to $\pi,F$ is given by
\begin{align*}
L(s,\pi\times F)=\sum_{m_1,m_2\geq 1}\frac{A_{\pi}(m_1,m_2)A_{\varpi}(m_1,m_2)}{(m_1^2m_2)^s},
\end{align*}
absolutely convergent for $\mathrm{Re}(s)>1$. 
With the above notation, we state the main theorem:
\begin{thm} \label{T1}
Let $q\geq 1$ be a large prime, and $\Omega\subseteq \mathfrak{a}_{\mathbb{C}}^{*}$.
Let $g$ be a holomorphic Hecke cusp form for $\mathrm{SL}_2(\mathbb{Z})$ of weight $k \equiv 0 \pmod{4}$.
Then we have
\begin{align*}
\sum_{\substack{F\subseteq L_{\mathrm{cusp}}^2(\Gamma_{0}(q)\textbackslash  \mathbb{H}^3)\\ \mu_{F}\in\Omega}}L(1/2,g \times F)L(1/2,\pi \times F)=q^{2}L(1,g \times \pi)+O\left(q^{\frac{11}{6}+\varepsilon}\right).
\end{align*}
\end{thm}

\begin{coro}
Under the notation of Theorem \ref{T1}, there exist infinitely many Hecke-Maa{\ss} cusp forms $F$ for $\mathrm{SL}_3(\mathbb{Z})$ such that
\[
L(1/2,g\times F)L(1/2,\pi \times F)\neq 0.
\]
\end{coro}

Throughout the paper, $\varepsilon$ denotes an arbitrarily small positive constant, which may vary from line to line. 
Moreover, following the standard convention in analytic number theory, we write $e(x)=e^{2\pi ix}$.

\section{Prerequisites}\label{secOutlin}
In this section, we collect some technical tools and background results that will be used in the proof of Theorem \ref{T1}.

\subsection{Automorphic forms and $L$-functions.}
We now recall the necessary properties of the automorphic forms and their associated $L$-functions.

Let $g$ be a holomorphic Hecke cusp form of weight $k$ for $\mathrm{SL}_2(\mathbb{Z})$. Its Fourier expansion is given by
\begin{align*}
g(z)=\sum^{\infty}_{n=1}\lambda_g(n)n^{\frac{k-1}{2}} e(nz)
\end{align*}
for $\mathrm{Im}(z)>0$, where $\lambda_g(n)\in\mathbb{R}$, and the Fourier coefficients satisfy the Hecke relation
\begin{align*}
\lambda_g(m)\lambda_g(n)=\sum_{d|(m,n)}\lambda_g\!\left(\frac{mn}{d^2}\right)
\end{align*}
with Deligne's bound $\lambda_g(n)\ll_{\varepsilon} n^{\varepsilon}$. The $L$-function associated to $g$ is defined by
\begin{align*}
L(s,g)=\sum_{n=1}^{\infty}\frac{\lambda_g(n)}{n^{s}}, \qquad \mathrm{Re}(s)>1,
\end{align*}
which admits analytic continuation to the complex plane and satisfies a functional equation of the form
\begin{align*}
\Lambda(s,g)=\varepsilon_g\Lambda(1-s,g),
\end{align*}
where $|\varepsilon_g|=1$ and
\begin{align*}
\Lambda(s,g)=(2\pi)^{-s}\Gamma\!\left(s+\frac{k-1}{2}\right)L(s,g).
\end{align*}

Let $\pi$ be a Hecke-Maa{\ss} cusp form on $\mathrm{GL}_3$ with Fourier coefficients $A_{\pi}(m,n)$ as introduced earlier. 
We define the normalized Hecke eigenvalues $\lambda_{\pi}(n)=A_{\pi}(1,n)$.
Known bounds toward the generalized Ramanujan conjecture imply that
\begin{align*}
A_{\varpi}(1,n)\ll n^{\theta_1+\varepsilon}
\end{align*}
for some $0\leq\theta_1\leq\frac{1}{2}$. 
The standard $L$-function attached to $\pi$ is
\begin{align*}
L(s,\pi)=\sum_{n=1}\frac{\lambda_{\pi}(n)}{n^s}
\end{align*}
for $\mathrm{Re}(s)>1$.
Let $\lambda_{\pi}=(\lambda_1,\lambda_2,\lambda_3)$ be the Langlands parameters of $\pi$,
satisfying $\lambda_1+\lambda_2+\lambda_3=0$. The completed $L$-function is defined by
\begin{align*}
\Lambda(s,\pi)=\pi^{-\frac{3s}{2}}\prod_{j=1}^3\Gamma\!\left(\frac{s+\lambda_j}{2}\right)L(s,\pi),
\end{align*}
which satisfies $\Lambda(s,\pi)=\varepsilon_{\pi}\Lambda(1-s,\pi)$ with $|\varepsilon_\pi|=1$.

To analyze the central values, 
we use approximate functional equations derived from the functional equations via contour shifting;
see \cite[Theorem 5.3]{IK}.
\begin{lemma}\label{L24}
Let $g$ be a Hecke cusp form for $\mathrm{SL}_2(\mathbb{Z})$, and let $F$ be a cuspidal automorphic representation on $\mathrm{GL}_3(\mathbb{Z})$ of level $q$. Then we have
\begin{align*}
&\hskip-0.5cmL(1/2,g\times F)=\\
&\hskip1cm\sum_{m \geq 1}\frac{\lambda_g(m)A_{\varpi}(1,m)}{m^{\frac{1}{2}}}\mathrm{V}\!\left(\frac{m}{\sqrt{\mathcal{Q}_1}}\right)+\varepsilon(g\times F)\sum_{n \geq 1}\frac{\lambda_g(n)\overline{A_{\varpi}(1,n)}}{n^{\frac{1}{2}}}\overline{\mathrm{V}}\!\left(\frac{n}{\sqrt{\mathcal{Q}_1}}\right),
\end{align*}
where the analytic conductor
\[\mathcal{Q}_1\asymp q^2,\] 
up to archimedean factors, and the weight functions are defined by
\[
\mathrm{V}(y)=\frac{1}{2\pi i}\int_{(3)}y^{-u}\frac{\gamma(g\times F,\frac{1}{2}+u)}{\gamma(g\times F,\frac{1}{2})}\frac{du}{u},\qquad
\overline{\mathrm{V}}(y)=\frac{1}{2\pi i}\int_{(3)}y^{-u}\frac{\gamma(g\times F,\frac{1}{2}-u)}{\gamma(g\times F,\frac{1}{2})}\frac{du}{u}.
\]
Here
\[
\gamma(g\times F,s)
=(2\pi)^{-3s}\prod_{i=1}^3
\Gamma\!\left(\frac{s+\frac{k-1}{2}+\mu_i}{2}\right),
\]
where $\mu_i$ denote the archimedean Satake parameters of $F$
satisfying $\mu_1+\mu_2+\mu_3=0$, and
\begin{align*}
\mathrm{V}^{(j)}(y),\overline{\mathrm{V}}^{(j)}(y)\ll_{A,j}(1+y)^{-A}
\end{align*}
for any $A>0$.
\end{lemma}
\begin{lemma}\label{L25}
Let $\pi$ and $F$ be cuspidal automorphic representations of $\mathrm{SL}_3(\mathbb{Z})$ of level $q$, we have
\begin{align*}
&L(1/2,\pi\times F)=\\
&\hskip1cm\sum_{m \geq 1}\frac{\lambda_{\pi}(m)A_{\varpi}(1,m)}{m^\frac{1}{2}}\mathrm{W}\!\left(\frac{m}{\sqrt{\mathcal{Q}_2}}\right)+
\varepsilon{(\pi\times F)}\sum_{n\geq 1}\frac{\lambda_{\pi}(n)\overline{A_{\varpi}(1,n)}}{n^\frac{1}{2}}\overline{\mathrm{W}}\!\left(\frac{n}{\sqrt{\mathcal{Q}_2}}\right),
\end{align*}
where the analytic conductor
\[\mathcal{Q}_2\asymp q^3,\] 
up to archimedean factors, and the weight functions are defined by
\[
\mathrm{W}(y)=\frac{1}{2\pi i}\int_{(3)}y^{-u}\frac{\gamma(\pi\times F,\frac{1}{2}+u)}{\gamma(\pi\times F,\frac{1}{2})}\frac{du}{u},\qquad
\overline{\mathrm{W}}(y)=\frac{1}{2\pi i}\int_{(3)}y^{-u}\frac{\gamma(\pi\times F,\frac{1}{2}-u)}{\gamma(\pi\times F,\frac{1}{2})}\frac{du}{u}.
\]
Here
\[
\gamma(\pi\times F,s)
=
\prod_{i=1}^3\prod_{j=1}^3
\pi^{s+\lambda_i+\mu_j}
\Gamma\!\left(\frac{s+\lambda_i+\mu_j}{2}\right),
\]
where $\lambda_i$ and $\mu_j$ denote the archimedean Satake parameters
of $\pi$ and $F$, respectively, and 
\[
\mathrm{W}^{(j)}(y),\overline{\mathrm{W}}^{(j)}(y)\ll_{A,j}(1+y)^{-A}
\]
for any $A>0$.
\end{lemma}

\subsection{The $\mathrm{GL}_3$ Kuznetsov Formula}
We give a recap of the Kuznetsov formula for $\mathrm{GL}_3$,
which generalizes the classical Bruggeman-Kuznetsov formula for $\mathrm{GL}_2$ \cite{K}. 
To state the formula, we require the following definitions.

Let $\{\varpi\}$ be an orthonormal basis of automorphic forms of level $q$,
cuspidal or Eisenstein series, containing all cuspidal newvectors.
We write $\int_{q}d\varpi$ or the corresponding spectral sum/integral.
The Fourier coefficients of $\varpi$ and its contragredient $\tilde{\varpi}$ are related by $A_{\varpi}(m,n)=A_{\tilde{\varpi}}(m,n)$ for any $(mn,q)=1$.
Then, we give the definition of an inner product on $(0, \infty)^2$,
\begin{align*}
\langle f_1,f_2 \rangle:=\int_{0}^{\infty}\int_{0}^{\infty}f_1(x,y)\overline{f_2(x,y)}\frac{dx_1dx_2}{(x_1x_2)^3}.
\end{align*}
For integers $D_1|D_2$, we define the $\mathrm{GL}_3$ Kloosterman sum
\[
\tilde{S}(m_1,n_1,n_2;D_1,D_2)=\sum_{\substack{C_1\bmod D_1 \\ C_2\bmod D_2 \\ (C_1,D_1)=(C_2,D_2/D_1)=1}}e{\left(n_{1}\frac{\tilde{C_{1}}C_{2}}{D_{1}}+n_{2}\frac{\tilde{C_{1}}}{D_{2}/D_{1}}+m_{1}\frac{C_{1}}{D_{1}}\right)}
\]
with the inverses taken with respect to the respective moduli. 
To incorporate the level structure, the modified Kloosterman sum is defined by
\begin{align}\label{21}
\notag&\hskip -15.2cm \tilde{S}^{(q)}(m_1,m_2,n_1,n_2;D_1,D_2)=\\
\hskip -0.2cm \sum_{\substack{B_1,C_1 \bmod D_1 \\ B_2,C_2 \bmod D_2 \\ B_1B_2+D_1C_2+D_2C_1\equiv 0\bmod D_1D_2 \\ (B_i,C_i,D_i)=1,i=1,2,q|B_1}}
e{\left(\frac{m_{1}B_{1}+n_{1}(Y_{1}D_{2}-Z_{1}B_{2})}{D_{1}}+\frac{m_{2}B_{2}+n_{2}(Y_{2}D_{1}-Z_{2}B_{1})}{D_{2}}\right)},
\end{align}
for $q|D_1,q|D_2$, where $Y_i,Z_i$ satisfy $Y_iD_i+Z_iC_i\equiv 1 \bmod D_i$ for $i=1,2$.
Note that the condition $q \mid B_1$ is redundant, since $q \mid B_1B_2$ already follows from $q \mid D_1$ and $q \mid D_2$.
With above notation, we now state the Kuznetsov formula for automorphic form $F$ with respect to the level aspect:
\begin{lemma}\label{L21}
Let $H:(0,\infty)^2\to \mathbb{C}$ be a smooth compactly supported test function satisfying the symmetry condition $H^{*}(y_1,y_2)=H(y_2,y_1)$. 
For $q,n_{1},n_{2},m_{1},m_{2} \in \mathbb{N}$, 
\[
\int_{(q)}\frac{\overline{A_{\varpi}(n_1,n_2)}\ A_{\varpi}(m_1,m_2)}{\mathcal{N}(\varpi)}\left|\left\langle \widetilde{W}_{\mu_{F}},H \right\rangle\right|^2d\varpi=\Delta+\Sigma_4+\Sigma_5+\Sigma_6,
\]
where
\begin{align*}
&\Delta=\delta_{n_1,m_1}\delta_{n_2,m_2}\|H\|^2,\\
&\Sigma_4=\sum_{\epsilon=\pm 1}\sum_{\substack{qD_2|D_1 \\ n_2D_1=m_1D^2_2}}\frac{\tilde{S}(\epsilon m_2,n_2,n_1,D_2,D_1)}{D_1D_2}\tilde{\mathcal{J}}_{\epsilon;H^*}{\left(\sqrt{\frac{n_1n_2m_2}{D_1D_2}}\right)},\\
&\Sigma_5=\sum_{\pm \epsilon}\sum_{\substack{q|D_1|D_2 \\ n_1D_2=m_2D^{2}_1}}\frac{\tilde{S}(\epsilon m_1,n_1,n_2,D_1,D_2)}{D_1D_2}\tilde{\mathcal{J}}_{\epsilon;H}{\left(\sqrt{\frac{n_1n_2m_1}{D_1D_2}}\right)},\\
&\Sigma_6=\sum_{\epsilon_1,\epsilon_2=\pm}\sum_{q|D_1,q|D_2}\frac{S^{(q)}(\epsilon_2m_2,\epsilon_1m_1,n_1,n_2,D_1,D_2)}{D_1D_2}\mathcal{J}_{\epsilon;H}{\left(\frac{\sqrt{n_2m_1D_1}}{D_2},\frac{\sqrt{n_1m_2D_2}}{D_1}\right)},
\end{align*}
where the integral transforms defined by
\begin{align*}
&\tilde{\mathcal{J}}_{\epsilon,H}(A) = A^{-2} \int_0^\infty \int_0^\infty \int_{-\infty}^\infty \int_{-\infty}^\infty e{(-\epsilon A x_1 y_1)} e{\left(y_2 \cdot \frac{x_1 x_2}{x_1^2 + 1}\right)}e{\left(\frac{A}{y_1 y_2} \cdot \frac{x_2}{x_1^2 + x_2^2 + 1}\right)} \\
&\hskip 1.7cm \times H{\left( y_2 \cdot \frac{\sqrt{x_1^2 + x_2^2 + 1}}{x_1^2 + 1}, \frac{A}{y_1 y_2} \cdot \frac{\sqrt{x_1^2 + 1}}{x_1^2 + x_2^2 + 1} \right)} \overline{H(A y_1, y_2)} \, dx_1 dx_2 \frac{dy_1 dy_2}{y_1 y_2^2}, \\
&\mathcal{J}_{\epsilon,H}(A_1, A_2) = (A_1 A_2)^{-2} \int_0^\infty \int_0^\infty \int_{-\infty}^\infty \int_{-\infty}^\infty \int_{-\infty}^\infty e{\left(-\epsilon_1 A_1 x_1 y_1 - \epsilon_2 A_2 x_2 y_2\right)} \\
&\hskip 1cm \times e{\left(-\frac{A_2}{y_2} \cdot \frac{x_1 x_3 + x_2}{x_3^2 + x_2^2 + 1}\right)} e{\left(-\frac{A_1}{y_1} \cdot \frac{x_2(x_1 x_2 - x_3) + x_1}{(x_1 x_2 - x_3)^2 + x_1^2 + 1}\right)} \overline{H(A_1 y_1, A_2 y_2)} \\ 
&\hskip 0.5cm \times H{\left( \frac{A_2}{y_2} \cdot \frac{\sqrt{(x_1 x_2 - x_3)^2 + x_1^2 + 1}}{x_3^2 + x_2^2 + 1}, \frac{A_1}{y_1} \cdot \frac{\sqrt{x_3^2 + x_2^2 + 1}}{(x_1 x_2 - x_3)^2 + x_1^2 + 1} \right)} dx_1 dx_2 dx_3 \frac{dy_1 dy_2}{y_1 y_2},
\end{align*}
and the Whittaker function as in \cite[Eqn. 2.15]{B} by its double Mellin transform
\begin{align*}
&\widetilde{W}_\mu(y_1, y_2) = \frac{y_1 y_2 \pi^{\frac{3}{2}}}{\left| \Gamma{\left( \frac{1}{2}(1 + i\Im(\mu_1 + 2\mu_2)) \right)} \Gamma{\left( \frac{1}{2}(1 + i\Im(\mu_1 - \mu_2)) \right)} \Gamma{\left( \frac{1}{2}(1 + i\Im(2\mu_1 + \mu_2)) \right)} \right|} \\ 
&\hskip 2cm \times \frac{1}{(2\pi i)^2} \int_{(1)} \int_{(1)} \frac{ \prod_{j=1}^3 \Gamma{\left( \frac{1}{2}(s_1 + \mu_j) \right)} \prod_{j=1}^3 \Gamma{\left( \frac{1}{2}(s_2 - \mu_j) \right) }}{4\pi^{s_1 + s_2} \Gamma{\left( \frac{1}{2}(s_1 + s_2) \right)}} y_1^{-s_1} y_2^{-s_2} \, ds_1 ds_2,
\end{align*}
where $A_1,A_2>0$, $\epsilon \in \{\pm 1\}$ or $\{\pm 1\}^2$ and $\delta_{m,n}$ is the Kronecker delta.
Moreover, the normalization factor (see \cite{B} for detailed definition) is
\begin{align}\label{Npi}
\mathcal{N}(\varpi)=[\mathrm{SL}_3(\mathbb{Z}):\Gamma_0(q)]\cdot \mathrm{Res}_{s=1} L(s,\varpi\times \tilde{\varpi}),
\end{align}
and it follows from \cite[Theorem 2]{Li10} that
\begin{align}\label{Npia}
\mathcal{N}(\varpi)\ll q^2(q(1+|\mu_{F}|))^{\varepsilon}.
\end{align}
\end{lemma}
\subsubsection{Kloosterman sum}
For specific Kloosterman sums defined as in \eqref{21}, we first present their individual bound
\begin{align*}
&S^{(q)}(m_1,m_2,n_1,n_2,D_1,D_2)\ll \\
&\hskip 3cm (D_1D_2)^{\frac{1}{2}+\varepsilon}\cdot((D_1,D_2)\cdot(m_1n_1,[D_1,D_2])\cdot(m_2n_2,[D_1,D_2]))^\frac{1}{2}  
\end{align*}
as shown in \cite[Eqn.(4.1)]{BBM}. If $q\| D_1$ and $q\| D_2$, and the Kloosterman sum satisfies
\[
S^{(q)}(m_1,m_2,n_1,n_2,D_1,D_2)=qS(\overline{q}m_1,\overline{q}m_2,n_1,n_2,D_1/q,D_2/q).
\]

To simplify notation, we abbreviated the notation $S^{(1)}$ to $S$ for the classical $\mathrm{GL}_3$ Kloosterman sum of level one.
In addition, two explicit evaluations are provided for special cases (which are essentially equivalent to \cite[Property 4.10]{DFG}).
\begin{lemma}\label{L221} We have
\begin{enumerate} 
\item[$(a)$]
$S^{(q)}(n_1, n_2, m_1, m_2; q, q^2) = r_q(m_1)S(n_2, m_2, q, q^2),$
\item[$(b)$]
$S^{(q)}(n_1, n_2, m_1, m_2; q^2, q) = r_q(m_2)S(n_1, m_1, q, q^2),$
\end{enumerate}
where $r_q(n)$ denote the Ramanujan sum.
\end{lemma}

In addition, as shown in the work of K{\i}ral and Nakasuji \cite{KJ},
the above Kloosterman sum admits a decomposition into a sum of products of two classical Kloosterman sums, namely
\begin{lemma}\label{L22}
Let $q \geq 1$, and let $S^{(q)}(m_1,m_2,n_1,n_2,D_1,D_2)$ be defined in \eqref{21}.
For any $D_1,D_2\geq 1$, we have
\begin{align}\label{2051}
&\hskip-0.7cmS(m_2,m_1,n_1,n_2,D_1,D_2)=\sum_{d|(D_1,D_2)}d\sum_{\substack{\gamma \bmod d \\ n_1D_2/d+m_1D_1\gamma/d\equiv 0 \bmod d}}\\
\notag&\hskip2.5cm S{\left(m_2,\frac{n_1D_2}{d^2}+\frac{m_1D_1\gamma}{d^2};\frac{D_1}{d}\right)}S{\left(n_2,\frac{n_1D_2\overline{\gamma}}{d^2}+\frac{m_1D_1}{d^2};\frac{D_2}{d}\right)}.
\end{align}
In addition, for $m_2, n_2, D_1, D_2\geq 1$ satisfying $(m_2n_2,q)=1$ and $(D_1D_2,q)=1$, one has
\[
S(m_1,m_2,n_1,n_2,qD_1,qD_2)=qS(n_1,\overline{q}m_2D_2;D_1)S(m_1,\overline{q}n_2D_1;D_2).
\]
\end{lemma}
\subsubsection{Weight functions in Kuznetsov formula}
For later convenience, we collect the properties of the integral transforms $\mathcal{J}_{\epsilon,H}$ and $\tilde{\mathcal{J}}_{\epsilon,H}$ (\cite[Lemma 3]{BBM}):
\begin{lemma}\label{L23}
We have
\begin{enumerate}
    \item[1.] \(\tilde{\mathcal{J}}_{\epsilon; H}(A) = 0 \) unless that \(A\gg 1\), in which case \(\tilde{\mathcal{J}}_{\epsilon,H}(A)\ll 1\), where the implied constant depends only on $H$.
    \item[2.] \(\mathcal{J}_{\epsilon,H}(A_1,A_2)=0\) unless that \(\min(A_1A_2^2,A_2A_1^2)\gg 1\). For $A_1,A_2$ in the range,
    \[
    \frac{\partial^i\partial^j}{\partial A_1^i \partial A_2^j}\mathcal{J}_{\epsilon,H}(A_1,A_2)\ll(A_1A_2)^\varepsilon(A_1^\frac{1}{3}A_2^\frac{2}{3})^i(A_1^\frac{2}{3}A_2^\frac{1}{3})^j.
    \] 
\end{enumerate}
for any $i,j \in \mathbb{N}_0$
\end{lemma}

\subsection{Auxiliary Lemma}
We record the following Wilton-type bound for cusp forms, as can be referenced in, for instance, \cite[Lemma 2.1.]{Y}.
\begin{lemma}\label{L27}
(Wilton-type bound) Let $h\in \mathbb{C}^{\infty}(\mathbb{R}^{\times,+})$ be a Schwartz class function that vanishes in a neighborhood of zero. 
For any newform $f\in \mathcal{B}_k^*$,
\[
\sum_{n\geq 1}\lambda_f(n)e(n\alpha)h{\left(\frac{n}{X}\right)}\ll_{\varepsilon} X^{\frac{1}{2}+\varepsilon}
\]
uniformly in any $\alpha\in \mathbb{R}$, where the implied constant depends on $\varepsilon$ and the parameter $k$ of $f$.
\end{lemma}

\section{Proof of Theorem \ref{T1}}\label{301}
In this section, we establish the proof of Theorem \ref{T1}. Let $F\subseteq L^{2}_{\mathrm{cusp}}(\Gamma_0(q)\textbackslash  \mathbb{H}^3)$ be a cuspidal automorphic representation.
Assume that $\varpi\in F$ is a newvector. Our main focus is
\begin{align}\label{301}
\mathcal{S}(q):=\sum_{\substack{F\subseteq L_{\mathrm{cusp}}^2(\Gamma_{0}(q)\textbackslash  \mathbb{H}^3)\\ \mu_{F}\in\Omega}}L(1/2,g\times F)L(1/2,\pi \times F)
\end{align}
Applying the approximate functional equations in Lemma \ref{L24} and Lemma \ref{L25},
we decompose \eqref{301} into
\begin{align}\label{31}
\mathcal{S}(q)=&\sum_{m,n\geq 1}\frac{\lambda_g(n)\lambda_{\pi}(m)}{\sqrt{mn}}\sum_{\substack{F\subseteq L_{\mathrm{cusp}}^2(\Gamma_{0}(q)\textbackslash  \mathbb{H}^3)\\ \mu_{F}\in\Omega}}A_{\varpi}(1,n) \overline{A_{\varpi}(1,m)}\mathrm{V}\!\left(\frac{n}{\sqrt{\mathcal{Q}_1}}\right)\overline{\mathrm{W}}\!\left(\frac{m}{\sqrt{\mathcal{Q}_2}}\right)\notag\\ 
&+\sum_{m,n\geq 1}\frac{\lambda_g(n)\lambda_{\pi}(m)}{\sqrt{mn}}\sum_{\substack{F\subseteq L_{\mathrm{cusp}}^2(\Gamma_{0}(q)\textbackslash  \mathbb{H}^3)\\ \mu_{F}\in\Omega}}\overline{A_{\varpi}(1,n)} A_{\varpi}(1,m)\overline{\mathrm{V}}\!\left(\frac{n}{\sqrt{\mathcal{Q}_1}}\right)\mathrm{W}\!\left(\frac{m}{\sqrt{\mathcal{Q}_2}}\right)\notag\\ 
&+\sum_{m,n\geq 1}\frac{\lambda_g(n)\lambda_{\pi}(m)}{\sqrt{mn}}\sum_{\substack{F\subseteq L_{\mathrm{cusp}}^2(\Gamma_{0}(q)\textbackslash  \mathbb{H}^3)\\ \mu_{F}\in\Omega}}\overline{A_{\varpi}(1,n)} \overline{A_{\varpi}(1,m)}\overline{\mathrm{V}}\!\left(\frac{n}{\sqrt{\mathcal{Q}_1}}\right)\overline{\mathrm{W}}\!\left(\frac{m}{\sqrt{\mathcal{Q}_2}}\right)\notag\\
&+\sum_{m,n\geq 1}\frac{\lambda_g(n)\lambda_{\pi}(m)}{\sqrt{mn}}\sum_{\substack{F\subseteq L_{\mathrm{cusp}}^2(\Gamma_{0}(q)\textbackslash  \mathbb{H}^3)\\ \mu_{F}\in\Omega}}A_{\varpi}(1,n) A_{\varpi}(1,m)\mathrm{V}\!\left(\frac{n}{\sqrt{\mathcal{Q}_1}}\right)\mathrm{W}\!\left(\frac{m}{\sqrt{\mathcal{Q}_2}}\right)\notag\\
=&\mathcal{S}_1+\mathcal{S}_2+\mathcal{S}_3+\mathcal{S}_4.
\end{align}
By the rapid decay of the weight functions, as stated in Lemma \ref{L24} and \ref{L25},
the sums are effectively restricted to
\begin{align}\label{re31}
m\ll q^{\frac{3}{2}},\quad n\ll q.    
\end{align}
Observe that the four terms are analogous; they differ only in the choice of weight functions,
which enjoy the identical decay properties. Hence it suffices to handle the term $\mathcal{S}_1$ in \eqref{31} in detail;
the treatment of $\mathcal{S}_2,\mathcal{S}_3,\mathcal{S}_4$ follows by identical arguments and will be omitted.
In particular, we focus on the following
\[
\sum_{m,n \geq 1}\frac{\lambda_g(n)\lambda_{\pi}(m)}{\sqrt{mn}}\sum_{\substack{F\subseteq L_{\mathrm{cusp}}^2(\Gamma_{0}(q)\textbackslash  \mathbb{H}^3)\\ \mu_{F}\in\Omega}}A_{\varpi}(1,n)\overline{A_{\varpi}(1,m)}\mathrm{V}\!\left(\frac{n}{\sqrt{\mathcal{Q}_1}}\right)\overline{\mathrm{W}}\!\left(\frac{m}{\sqrt{\mathcal{Q}_2}}\right).
\]
We apply the spectral averaging procedure for the $\mathrm{GL}_3$ Kuznetsov formula, following \cite[Section 5.1]{BBM}, together with the estimate \eqref{Npia} for \eqref{Npi}. This yields
\begin{align*}
&\ll q^{2+\varepsilon}\sum_{m,n \geq 1}\frac{\lambda_g(n)\lambda_{\pi}(m)}{\sqrt{mn}}\max_{i}\max_{|t|\ll q^{\varepsilon}}
\sum_{\substack{F\subseteq L_{\mathrm{cusp}}^2(\Gamma_{0}(q)\textbackslash  \mathbb{H}^3)\\ \mu_{F}\in\Omega}}|\langle\tilde{W}_{\mu_{F}},H_i\rangle|^2\\
&\hskip5cm\times\frac{A_{\varpi}(1,n)\overline{A_{\varpi}(1,m)}}{\mathcal{N}(\varpi)}\mathrm{V}\!\left(\frac{n}{\sqrt{\mathcal{Q}_1}}\right)\overline{\mathrm{W}}\!\left(\frac{m}{\sqrt{\mathcal{Q}_2}}\right),
\end{align*}
where $\{H_1,H_2,\dots,H_s\}$ is a finite collection of smooth compactly supported functions, such that $\sum_i|\langle\tilde{W}_{\mu_{F}},H_i\rangle|\gg 1$ for $\mu_F\in\Omega$, as shown in \cite[Lemma 1]{BBM}.
By re-arranging, the right-hand side is thus
\begin{align*}
&q^{2+\varepsilon}\sum_{m,n \geq 1}\frac{\lambda_g(n)\lambda_{\pi}(m)}{\sqrt{mn}}\mathrm{V}\!\left(\frac{n}{\sqrt{\mathcal{Q}_1}}\right)\overline{\mathrm{W}}\!\left(\frac{m}{\sqrt{\mathcal{Q}_2}}\right)
\\&\hskip3.5cm\times\max_{i}\int_{(q)}\frac{A_{\varpi}(1,n),\overline{A_{\varpi}(1,m)}}{\mathcal{N}(\varpi)}|\langle\tilde{W}_{\mu_{F}},H_i\rangle|^2d F.
\end{align*}
Denote the above spectral sum expression by $\Xi$, and by employing the $\mathrm{GL}_3$ Kuznetsov formula in Lemma \ref{L21}, we get
\[
\Xi=\Delta_{m,n}+\Sigma_4+\Sigma_5+\Sigma_6,
\]
where 
\[
\Delta_{m,n}=q^{2+\varepsilon}\max_{i}\|H_i\|^2\sum_{m,n\geq 1}\frac{\lambda_g(n)\lambda_{\pi}(m)}{\sqrt{mn}}\mathrm{V}\!\left(\frac{n}{\sqrt{\mathcal{Q}_1}}\right)\overline{\mathrm{W}}\!\left(\frac{m}{\sqrt{\mathcal{Q}_2}}\right)\mathbf{1}_{m=n}
\]
is the diagonal term of the decomposition, and $\Sigma_4,\Sigma_5,\Sigma_6$ are off-diagonal sums.
Explicitly, these terms are given by
\begin{align*}
&\Sigma_4=q^{2+\varepsilon}\max_{i}\sum_{\epsilon=\pm 1}\sum_{m,n\geq 1}\frac{\lambda_g(n)\lambda_{\pi}(m)}{\sqrt{mn}}\sum_{\substack{qD_2|D_1 \\ mD_1=D_2^2}}\frac{\tilde{S}(\epsilon n,m,1,D_2,D_1)}{D_1D_2}\\
&\hskip6cm\times\mathcal{J}_{\epsilon,H}{\left(\sqrt{\frac{mn}{D_1D_2}}\right)}\mathrm{V}\!\left(\frac{n}{\sqrt{\mathcal{Q}_1}}\right)\overline{\mathrm{W}}\!\left(\frac{m}{\sqrt{\mathcal{Q}_2}}\right), \\
&\Sigma_5=q^{2+\varepsilon}\max_{i}\sum_{\epsilon=\pm 1}\sum_{m,n\geq 1}\frac{\lambda_g(n)\lambda_{\pi}(m)}{\sqrt{mn}}\sum_{\substack{q|D_1|D_2 \\ D_2=nD_1^2}}\frac{\tilde{S}(\epsilon ,1,m,D_1,D_2)}{D_1D_2}\\
&\hskip6cm\times\mathcal{J}_{\epsilon,H}{\left(\sqrt{\frac{m}{D_1D_2}}\right)}\mathrm{V}\!\left(\frac{n}{\sqrt{\mathcal{Q}_1}}\right)\overline{\mathrm{W}}\!\left(\frac{m}{\sqrt{\mathcal{Q}_2}}\right), \\
\end{align*}
and
\begin{align}\label{3111}
&\Sigma_6=q^{2+\varepsilon}\max_{i}\sum_{\epsilon_1,\epsilon_2=\pm 1}\sum_{m,n\geq 1}\frac{\lambda_g(n)\lambda_{\pi}(m)}{\sqrt{mn}}\sum_{q|D_1,q|D_2}\frac{S^{(q)}(\epsilon_2 n,\epsilon_1 ,1,m,D_1,D_2)}{D_1D_2}\\ 
\notag&\hskip6cm\times\mathcal{J}_{\epsilon,H}{\left(\frac{\sqrt{mD_1}}{D_2},\frac{\sqrt{nD_2}}{D_1}\right)}\mathrm{V}\!\left(\frac{n}{\sqrt{\mathcal{Q}_1}}\right)\overline{\mathrm{W}}\!\left(\frac{m}{\sqrt{\mathcal{Q}_2}}\right).
\end{align}

\subsection{Treatment of $\Delta_{m,n}$}\label{DeltaT}
We first consider the diagonal term $\Delta_{m,n}$.
Setting $m=n$, we obtain
\begin{align*}
\Delta_{m,m}=q^{2+\varepsilon}\sum_{m\geq 1}\frac{\lambda_g(m)\lambda_{\pi}(m)}{m}\mathrm{V}\!\left(\frac{n}{\sqrt{\mathcal{Q}_1}}\right)\overline{\mathrm{W}}\!\left(\frac{m}{\sqrt{\mathcal{Q}_2}}\right).
\end{align*}
Here, we insert the Mellin representations of the weight functions given in Lemma \ref{L24} and Lemma \ref{L25}, together with Rankin-Selberg convolution, this yields
\begin{align*}
\Delta_{m,m}=q^{2+\varepsilon}\frac{1}{(2\pi i)^2}\int_{(3)}\int_{(3)}L(1+u+v,g\times \pi)
\mathcal{Q}_1^{\frac{u}{2}}\mathcal{Q}_2^{\frac{v}{2}}G_1(u)G_2(v)\frac{dudv}{uv},
\end{align*}
where
\begin{align*}
G_1(u)=\frac{\gamma(g\times F,\frac{1}{2}+u)}{\gamma(g\times F,\frac{1}{2})},\qquad G_2(v)=\frac{\gamma(\pi\times F,\frac{1}{2}-v)}{\gamma(\pi\times F,\frac{1}{2})}.
\end{align*}
By Stirling's formula in vertical strips, the ratios $G_1(u)$ and $G_2(v)$ are of polynomial growth in $|u|$ and $|v|$, and the integrand is absolutely convergent on $\mathrm{Re}(u)=\mathrm{Re}(v)=3$.

We now shift the contour in $u$ first, keeping $v$ fixed. Moving the $u$-contour from $\mathrm{Re}(u)=3$ to $\mathrm{Re}(u)=-\frac{1}{2}$, we encounter a simple pole at $u=0$ arising from the factor $1/u$. 
The residue at $u=0$ contributes 
\begin{align*}
q^{2+\varepsilon} \frac{1}{2\pi i}\int_{(3)}L(1+v, g \times \pi)\mathcal{Q}_2^{v/2}G_2(v)\frac{dv}{v}.
\end{align*}
No additional poles are encountered in the shift, 
since $L(s, g \times \pi)$ is holomorphic in $\mathrm{Re}(s) > 0$ except for a simple pole at $s=1$.
We next shift the contour in $v$ from $\mathrm{Re}(v)=3$ to $\mathrm{Re}(v)=-\frac{1}{2}$. 
During this shift, we encounter a simple pole at $v=0$, again from the factor $1/v$. The residue at $v=0$ gives the main term
\begin{align*}
q^{2} L(1, g \times \pi).
\end{align*}
It remains to bound the contributions from the shifted contours. 
On the lines $\Re(u)=-1/2$ and $\Re(v)=-1/2$, we use the convexity bound $L(1+u+v, g \times \pi) \ll q^\varepsilon$, together with 
\begin{align*}
\mathcal{Q}_1^{\frac{u}{2}} \ll q^{-\frac{1}{2}}, \qquad \mathcal{Q}_2^{\frac{v}{2}} \ll q^{-\frac{3}{4}},
\end{align*}
and the estimates with the polynomial growth of $G_1(u)$ and $G_2(v)$, we find that the contribution from the shifted contours is
\begin{align*}
\ll q^{-\frac{5}{4}+\varepsilon}(1+|u|+|v|)^{A}
\end{align*}
for some fixed $A>0$. It follows that the contribution from the shifted contours is $\ll q^{-5/4+\varepsilon}$, and hence
\begin{align*}
\Delta_{m,n}=q^2L(1,g\times\pi)+O(q^{\frac{3}{4}+\varepsilon}).
\end{align*}

\subsection{Treatment of $\Sigma_4,\Sigma_5$}\label{SigmaT}
For $\Sigma_4$, we simply observe that the conditions $qD_2|D_1$ and $mD_1=D_2^2$,
we may write $D_1=aqD_2$ for some $a\geq 1$, so that $D_2=amq$ and $D_1=ma^2q^2$.
Therefore $D_1D_2\geq m^2q^3$. On the other hand, it follows from \eqref{re31} that $mn\ll qm$,
so that (for sufficiently small $\varepsilon>0$ and sufficiently large $q$) we have $\Sigma_4$ is negligible by Lemma \ref{L24} (1).
For $\Sigma_5$, the conditions $q|D_1|D_2,D_2=nD_1^2$ implies $D_1D_2\geq q^3$, but $n\ll q$, which we have $\Sigma_5$ is negligible in similar argument.

\subsection{Treatment of $\Sigma_6$} In the next step, we will estimate the term $\Sigma_6$, one observes that the expression in \eqref{3111} is dominated by
\begin{align*}
&\Sigma_6\ll q^{2+\varepsilon}\sum_{\epsilon\in\{\pm 1\}^2}\sum_{m,n \geq 1}\frac{\lambda_g(n)\lambda_{\pi}(m)}{\sqrt{nm}}\sum_{q|D_1,q|D_2}\frac{S^{(q)}(n,1,1,m,D_1,D_2)}{D_1D_2}\notag\\
&\hskip5cm\times\mathcal{J}_{\epsilon,H}{\left(\frac{\sqrt{mD_1}}{D_2},\frac{\sqrt{nD_2}}{D_1}\right)}\mathrm{V}\!\left(\frac{n}{\sqrt{\mathcal{Q}_1}}\right)\overline{\mathrm{W}}\!\left(\frac{m}{\sqrt{\mathcal{Q}}_2}\right).
\end{align*}
Here, we proceed by examining the Kloosterman sum in the above sum. We analyze the dependence on the variable $m$, with respect to $q$.
It suffices to restrict to the non-degenerate case $(q,m)=1$;
the degenerate case where $q^{\eta}|m$ (for any $\eta\geq 1$) can be treated similarly and contributes less.
We may thus assume that $(q,m)=1$, in which the following subcases arise 
\begin{itemize}
\item[(a)] $q\|D_1, q\|D_2$.
\item[(b)] $q^2|D_1$ or $ q^2\|D_2$.
\end{itemize}
In Case $(b)$, by invoking Lemma \ref{L221}, which implies that the contribution of the term is zero.
In Case $(a)$, we write $D_1=qc_1c^*$ and $D_2=qc_2c^*$, with $c^*=(D_1/q,D_2/q)$ satisfies $(c_1,c_2)=1$.
Via \eqref{2051}, we deduce that
\begin{equation}\label{crefkfs}
\begin{aligned} 
&{S} ^{(q)}(n,1,1 ,m, q c_1 c^\ast,q c_2 c^\ast )=\\
&\hskip1cm q c^\ast  \sum_{\substack{\alpha\bmod {c^\ast }\\ c_2+c_1\alpha \equiv 0 \bmod {c^\ast} }} 
S\left(n,\frac{\overline{q}(c_2+c_1\alpha)}{c^\ast};c_1\right)
S\left(m,\frac{\overline{q}(c_2\overline{\alpha}+c_1)}{c^\ast};c_2\right).
\end{aligned}\end{equation}
Here, $\overline{\alpha}$ denotes the modular inverse of $\alpha$ modulo $c^*$, satisfying $\alpha \overline{\alpha} \equiv1 \bmod {c^\ast}$.
To continue the derivation, given the preceding summation constraints, one observes that $c^*$ is coprime with one of the two elements $\{c_1,c_2\}$.
Without loss of generality, we assume that $(c^*,c_1)=1$,
the summation condition on $\alpha$ implies that the expression in \eqref{crefkfs} simplifies to
\[ 
qc^\ast S(n,\overline{qc^\ast}c_2;c_1)S(m,\overline{qc^\ast}c_1;c_2),
\]
where if $c^*=1$, instead of the formula in \eqref{crefkfs}, we have
\[
{S}^{(q)}(n,1,1 ,m, q c_1 c^\ast,q c_2 c^\ast )=qS(n,\overline{q}c_2;c_1)
S(m,\overline{q}c_1;c_2).
\]
We first restrict attention to the dominant case where $c^*\ll q^{\varepsilon}$;
the complementary case $c^*\gg q^{\varepsilon}$ can be treated similarly and yields a smaller contribution, which we omit for brevity.
Then we replace $D_1,D_2$ by $c_1q,c_2q$, respectively. Let $\Phi(x,y;z_1,z_2)$ define as
\begin{align*}
\Phi_{\epsilon,H}(x,y;z_1,z_2)=\mathcal{J}_{\epsilon,H}{\left(\frac{\sqrt{xz_1}}{z_2\sqrt{q}},\frac{\sqrt{yz_2}}{z_1\sqrt{q}}\right)}\mathrm{V}\!\left(\frac{y}{\sqrt{\mathcal{Q}_1}}\right)\overline{\mathrm{W}}\!\left(\frac{x}{\sqrt{\mathcal{Q}_2}}\right),
\end{align*}
upon combining the properties of Lemma \ref{L23} (2), Lemma \ref{L24} and Lemma \ref{L25}, which satisfies that
\begin{align}\label{35}
\frac{\partial^i\partial^j}{\partial c_1^i\partial c_2^j}\Phi_{\epsilon,H}(m,n;c_1,c_2)\ll q^{\frac{i}{2}+\frac{j}{2}+\varepsilon}
\end{align}
for any $i,j>0$.
This reduces the interest to the evaluation of
\begin{align*}
\Sigma_6\ll q^{1+\varepsilon}\sum_{m,n\geq 1}\frac{\lambda_g(n)\lambda_{\pi}(m)}{\sqrt{mn}} \sum_{\substack{c_1c_2 \geq 1}}
\frac{S(m,\overline{q}c_1;c_2)S(n,\overline{q}c_2;c_1)}{c_1c_2}\Phi_{\epsilon,H}(m,n;c_1,c_2).
\end{align*}
By Lemma \ref{L23} (2), together with \eqref{re31}, the kernel $\mathcal{J}_{\epsilon,H}$ is supported only when
\begin{equation}\label{re32}\begin{gathered}
c_1\ll n^\frac{2}{3}m^\frac{1}{3}q^\varepsilon/q\ll q^{\frac{1}{6}+\varepsilon}, \\ 
c_2 \ll n^\frac{1}{3}m^\frac{2}{3}q^\varepsilon/q\ll q^{\frac{1}{3}+\varepsilon}
\end{gathered}\end{equation} 
Next, by applying Cauchy-Schwarz inequality to the core sum to get rid of the Fourier coefficients $\lambda_{\pi}(m)$, thereby yielding
\begin{align}\label{32}
\Sigma_6 \ll q^{\frac{1}{2}+\varepsilon}\cdot \mathcal{A}^{\frac{1}{2}}(q,\Phi_{\epsilon,H}),
\end{align}
where the multiple sum $\mathcal{A}(q,\Phi_{\epsilon,H})$ is defined by
\begin{equation}\label{36}\begin{aligned}
\mathcal{A}(q,\Phi_{\epsilon,H})=&\sum_{m\geq 1}\sum_{n_1,n_2\geq 1}\frac{\lambda_g(n_1)\lambda_g(n_2)}{\sqrt{n_1n_2}}
\sum_{\substack{c_1,c_1'\geq 1 \\ c_2,c_2'\geq 1}}S(n_1,\overline{q}c_2;c_1)\overline{S(n_2,\overline{q}c_2';c_1')}\\
&\times S(m,\overline{q}c_1;c_2)\overline{S(m,\overline{q}c_1';c_2')}
\Phi_{\epsilon_{\epsilon,H}}\left(\frac{m}{\mathcal{L}},\frac{n_1}{\mathcal{L}};c_1,c_2\right)\overline{\Phi_{\epsilon,H}\left(\frac{m}{\mathcal{L}},\frac{n_2}{\mathcal{L}};c_1',c_2'\right)}.
\end{aligned}\end{equation}
Here, $\mathcal{L}\geq 1$ is a parameter to be chosen very soon,
and we have utilized the known bound of \textit{Ramanujan on average} that $\sum_{n\leq x}A_{\varpi}(1,n)\ll x^{1+\varepsilon}$ for any $x>2$,
which directly follows from the Rankin-Selberg theory; see \cite[Remark 12.1.8]{Go} for detail.

Next, by applying the Poisson summation formula to the variable $\mathrm{m}$, with modulus $c_2c_2'$,
it follows that the $\mathrm{m}$-sum in \eqref{36} is essentially converted into a shape of
\begin{align*}
\frac{q^{\frac{3}{2}+\varepsilon}\mathcal{L}}{c_2c_2'}\sum_{|h|\ll\frac{(c_2c_2')^{1+\varepsilon}}{\mathcal{L}}}
\sum_{\alpha\bmod c_2c_2'}S(\alpha,\overline{q}c_1;c_2)\overline{S(\alpha,\overline{q}c_1';c_2')}e\!\left(-\frac{\alpha h}{c_2c_2'}\right).
\end{align*}
We choose $\mathcal{L}=q^{100\varepsilon}$, so that the contribution of \textit{non-zero frequencies} is negligible.
Indeed, under this choice one has
\[c_2c_2'\ll q^{\frac{3}{2}-\varepsilon}\mathcal{L},\]
for any $\varepsilon>0$.
Here, it suffices to consider the diagonal case, namely, $c_2=c_2'$, need to be taken into account (for the main contribution);
specifically, we would invoke the Chinese Remainder Theorem to decompose the exponential sum $\alpha$, with modulus $c_2\neq c_2'$,
which leads to sums of the form
\begin{align*}
\frac{q^{\frac{3}{2}+\varepsilon}\mathcal{L}}{c_2c_2'}\sum_{|h|\ll\frac{(c_2c_2')^{1+\varepsilon}}{\mathcal{L}}}
\sum_{\alpha_1\bmod c_2}S(\alpha_1,\overline{q}c_1;c_2)e\biggl(-\frac{\alpha_1\overline{c_2'}h}{c_2}\biggr)
\sum_{\alpha_2\bmod c_2'}\overline{S(\alpha_2,\overline{q}c_1';c_2')}e\biggl(-\frac{\alpha_2\overline{c_2}h}{c_2'}\biggr)
\end{align*}
which \textit{vanishes}. 
It follows that the total contribution of the off-diagonal terms $c_2\neq c_2'$ is negligible compared to that of the diagonal term.
We are thus left with the diagonal case $c_2=c_2'$.
In this setting, an entirely analogous argument with the help of Poisson to the variable $\mathrm{m}$, with modulus $c_1c_1'$, on the other hand, 
implies that the display \eqref{36} is also negligibly small, which yields the necessary condition that
\[
c_1\equiv c_1' \bmod c_2.
\]
Taking into account the size constraints $c_1, c_1' \ll c_2$, the congruence condition forces $c_1 = c_1'$. Consequently, we can infer that $\mathcal{A}(q,\Phi_{\epsilon,H})$ in \eqref{36} is transformed into
\begin{align*}
&\sum_{m\geq 1}\sum_{n_1,n_2\geq 1}\frac{\lambda_g(n_1)\lambda_g(n_2)}{\sqrt{n_1n_2}}\sum_{\substack{c_1,c_2\geq 1}}S(n_1,\overline{q}c_2;c_1)\overline{S(n_2,\overline{q}c_2;c_1)}\\
&\hskip3cm\times |S(m,\overline{q}c_1;c_2)|^2\Phi_{\epsilon,H}(m,n_1;c_1,c_2)\overline{\Phi_{\epsilon,H}(m,n_2;c_1,c_2)}.
\end{align*}
By swapping the order of summation, we can bound the expression in \eqref{36} more efficiently, resulting in an estimate of the following form
\begin{align*}
&\sum_{m\geq 1}\sum_{\substack{c_1, c_2\geq 1}}|S(m,\overline{q}c_1;c_2)|^2 \sum_{n_1,n_2\leq q}\frac{\lambda_g(n_1)\lambda_g(n_2)}{\sqrt{n_1n_2}}\\
&\hskip2.5cm\times S(n_1,\overline{q}c_2;c_1)\overline{S(n_2,\overline{q}c_2;c_1)}
\Phi(m,n_1;c_1,c_2)\overline{\Phi(m,n_2;c_1,c_2)}.
\end{align*}
Upon opening the Kloosterman sum and appealing to Lemma \ref{L27}, and the condition \eqref{re31} and \eqref{re32} on the above sum, with the estimate \eqref{35} for $\Phi$, we arrive at the display \eqref{36} is
\begin{align}\label{42}
\ll \sum_{m\ll q^{\frac{3}{2}+\varepsilon}}\sum_{\substack{c_1\ll q^{\frac{1}{6}+\varepsilon} \\ c_2\ll q^{\frac{1}{3}+\varepsilon}}}c_2^{1+\varepsilon}c_1^{2+\varepsilon}\ll q^{2+\frac{2}{3}+\varepsilon}\ll q^{\frac{8}{3}+\varepsilon}.
\end{align}
By combining the key relations from the equation in \eqref{36} and the estimate \eqref{42}, we are able to derive the following result
\[
\Sigma_6\ll q^{\frac{1}{2}+\varepsilon}{\left(q^{\frac{8}{3}+\varepsilon}\right)}^{\frac{1}{2}}\ll q^{\frac{11}{6}+\varepsilon}.
\]
Here, collecting the estimates from the contributions of the diagonal term $\Delta_{m,n}$ in section \ref{DeltaT} as well as the off-diagonal sums $\Sigma_4,\Sigma_5$ in section \ref{SigmaT} and $\Sigma_6$ above,
this completes the proof of Theorem \ref{T1}.

\section{Acknowledgments}
The author is grateful to Prof. Hou for help during the preparation of this paper.

\end{document}